\documentclass[a4paper,12pt]{article}
\usepackage{graphicx}
\usepackage{amssymb}
\usepackage{amsmath}
\input xy
\xyoption{all}
\begin{document}

\newcommand{\thm}[1][\!\!]{\textbf{Theorem {#1}.\;}}
\newcommand{\lem}[1][\!\!]{\textbf{Lemma {#1}.\;}}
\newcommand{\pf}[1][\!\!]{\textit{Proof {#1}.\;}}
\newcommand{\rk}[1][\!\!]{\textit{Remark {#1}.\;}}
\newcommand{\cor}[1][\!\!]{\textit{Corollary {#1}.\;}}
\newcommand{\defn}[1][\!\!]{\textit{Definition {#1}.\;}}
\newcommand{\qed}{$\square$}
\newcommand{\Nproj}{{\vphantom{N\;}}_{N\,}{\!\mathrm{proj}}}
\newcommand{\projN}{\mathrm{proj}_N}
\newcommand{\Stab}{\mathrm{Stab}}
\newcommand{\SL}{\mathrm{SL}}
\newcommand{\B}{\mathbf{B}}


\begin{titlepage}
\title{An equivariant covering map from the upper half plane to the complex plane minus a lattice}
\author{M. Batchelor$^a$ \and P. Brownlee\footnote{e-mail \texttt{pollybrownlee@gmail.com}} \and W. L. Woods\footnote{e-mail \texttt{wlw21@cam.ac.uk}}}

\maketitle

\begin{center}
$^a$ Corresponding Author\\
Dept. of Pure Mathematics and Mathematical Statistics,\\
Centre for Mathematical Sciences, University of Cambridge,\\
Wilberforce Road, Cambridge, CB3 0WB, United Kingdom\\
Tel:  +44 (0) 1223 765896, Fax: +44 (0) 1223 337920, Email: \texttt{mb139@cam.ac.uk}
\end{center}


\begin{abstract}

This paper studies a covering map $\varphi$ from the upper half plane to the complex plane with a triangular lattice excised.  This map is interesting as it factorises Klein's $J$ invariant. Its derivative has properties which are a slight generalisation of modular functions, and $(\varphi')^6$ is a modular function of weight 12. There is a homomorphism from the modular group $\Gamma$ to the affine transformations of the complex plane which preserve the excised lattice.  With respect to this action $\varphi$ is a map of $\Gamma$--sets. Identification of the excised lattice with the root lattice of $\mathfrak{sl}_3(\mathbb{C})$ allows functions familiar from the study of modular functions to be expressed in terms of standard constructions on representations of $\mathfrak{sl}_3(\mathbb{C})$.  
\end{abstract}
\end{titlepage}
\medskip

\section*{Introduction}

The work described here was inspired by the striking similarity of the following familiar spaces:
\begin{center}\includegraphics[height=48mm]{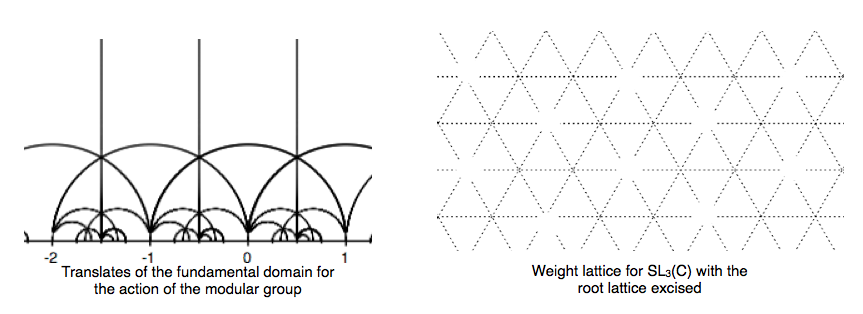}\end{center}
The set of paths on these two one-dimensional spaces are the same. Since the one on the left is a tree, and hence contractible, it must be the universal cover of the one on the right.

This paper studies the extension of this covering map of the one dimensional spaces to a map $\varphi$ of the upper half plane $\mathbb{H}$, to $\B$, the complex plane with points in a lattice excised:
\[
\B = \mathbb{C} \setminus \{m(\omega+1) + n(\omega^2-1) \mid m, n \in \mathbb{Z} \}
\]
where $\omega = \frac{1}{2}(1+i\sqrt{3})$.  The action of $\Gamma$ on this space is achieved via a homomorphism onto a subgroup $\overline{\Gamma}$ of the group of affine transformations of $\mathbb{C}$ which preserve $\B$.  

The uniformisation theorem guarantees the existence of a covering map \mbox{$\varphi : \mathbb{H} \to \B$}.  The principal result of this paper places this map in the context of well known functions:
\[
(\wp'(\varphi(z)))^2 = J(z)
\]
where $\wp$ is the Weierstrass $\wp$ function and $J$ is the Klein invariant.  Moreover, this map is equivariant, in the sense that the following diagram commutes.
\begin{center}
$\xymatrix{\Gamma\times\mathbb{H}\ar[r]^-{\text{action}} \ar[d]_{ \psi \times\varphi} & \mathbb{H}\ar[d]^\varphi \\
\overline{\Gamma}\times \B\ar[r]_-{\text{action}} & \B \\
}$
\end{center}
An immediate corollary of this is that $\varphi'$ has properties which are a slight generalisation of those required of modular functions:
\[
\varphi'(gz) = \ell(g)(cz+d)^2\varphi'(z),
\]
where $\ell$ is a homomorphism from $\Gamma$ to the sixth roots of unity in $\mathbb{C}$.

A tractable explicit presentation for $\varphi$ has proved elusive, but in terms of the group theory, $\B$ can be readily identified with a double coset space. $\mathrm{SL}_2(\mathbb{R})$ acts transitively on $\mathbb{H}$, so as usual, choice of a point $p \in \mathbb{H}$ identifies $\mathbb{H}$ with $\mathrm{SL}_2(\mathbb{R})/\Stab(p)$, where this stabiliser is taken in $\mathrm{SL}_2(\mathbb{R})$.  Let $N$ denote the kernel of the map $\mathrm{SL}_2(\mathbb{Z})\to \Gamma\to\overline{\Gamma}$.  Then

\[
\B \simeq N\backslash SL_2(\mathbb{R}) / \Stab(p).
\]

The map $\varphi$ implies a connection between modular functions and the representation theory of $\mathfrak{sl}_3(\mathbb{C})$.  This connection arises through recognising $\B$ as the real span of the weight space of $\mathfrak{sl}_3(\mathbb{C})$.  This allows familiar functions such as $\wp$ and modular $\lambda$ to be interpreted as operators arising from representations of $\mathfrak{sl}_3(\mathbb{C})$.

The paper is organised as follows.

\begin{enumerate}
	\item The group homomorphism $\psi : \Gamma \to \overline{\Gamma}$.
	\item Definition of $\varphi$.
	\item Group--theoretic description of $\varphi$.	
	\item The interpretation of  $\wp$ and modular $\lambda$ in terms of the representation theory of $\mathfrak{sl}_3(\mathbb{C})$.
	
\end{enumerate}

We would like to thank James Bridgwater for support during the period of this research.

\section{The group homomorphism $\psi$ from $\Gamma$ to the group of affine transformations of $\mathbb{C}$}

The modular group $\Gamma = \mathrm{SL}_2(\mathbb{Z})/\{\pm 1\}$ has a standard presentation with generators
\begin{center}
$S = \begin{pmatrix} 0&-1\\1&0 \end{pmatrix}, T = \begin{pmatrix} 1&1\\0&1 \end{pmatrix}$
\end{center}
which satisfy the relations $S^2 = (ST)^3 = 1$ only. Under the action of the modular group on $\mathbb{H}$,
\begin{center}
$\begin{pmatrix} a&b\\c&d \end{pmatrix} \cdot z = \dfrac{az + b}{cz + d}$,
\end{center}
$S$ stabilises $i$, $ST$ stabilises $\omega$, and $TS$ stabilises $\omega ^2$.

The group of affine transformations $\mathrm{Aff}(\mathbb{C})$ can be represented as the subgroup of matrices
\[
\mathrm{Aff}(\mathbb{C}) = \left\{ \begin{pmatrix} a&b\\0&1 \end{pmatrix} \right\}\subset \mathrm{GL}_2(\mathbb{C}),
\]
where the action on $\mathbb{C}$ is given by
\[
A(u) = \mathrm{proj}_1\left(A\begin{pmatrix} u\\1 \end{pmatrix}\right),
\]
where $\mathrm{proj}_1$ is projection onto the first component. Define a representation $\psi: \Gamma\to \mathrm{Aff}(\mathbb{C})\subset \mathrm{GL}_2(\mathbb{C})$ by

\begin{center}
$S\mapsto \overline{S} = \begin{pmatrix} -1&\omega + \omega^2 \\ 0&1 \end{pmatrix}$,
$T\mapsto \overline{T} = \begin{pmatrix}\omega&0 \\ 0&1 \end{pmatrix}$,
\end{center}
and let $\overline{\Gamma}$ be the image of $\Gamma$ under this representation. Note that $\overline{S}$ and $\overline{T}$ preserve the lattice $L = \mathbb{Z} \oplus \mathbb{Z}\omega$, thus $\overline{\Gamma}$ is a subgroup of $\mathrm{Aut}(L)$. ($\overline{\Gamma}$ is the group that will give a $\Gamma$-action on $\B$.)

\lem $\overline{\Gamma}$ has presentation $\langle \overline{S}, \overline{T} | \overline{S}^2, (\overline{S} \overline{T})^3, \overline{T}^6 \rangle$.

\pf $\overline{\Gamma}$ is actually a subgroup of $\mathrm{Aut}^+(L)$, the group of orientation-preserving automorphisms of $L$, as both $\overline{S}$ and $\overline{T}$ are rotations.

Start by finding a presentation for $\mathrm{Aut}^+(L)$. Any element $\sigma\in \mathrm{Aut}^+(L)$ is completely determined by its action on two points, so it suffices to track the pair $P = (0, 1)$. Note that $\mathrm{Aut}^+(L)$ is generated by three elements: a translation $a: z\mapsto z+1$, a translation $b: z\mapsto z+\omega$ and a rotation $c: z\mapsto \omega z$. By observing $P$ and its image $\sigma P$, we note the following six relations:

\begin{center}
$ba = ab, \quad ca = bc,\quad cb = a^{-1}bc,\quad ca^{-1} = b^{-1}c,\quad cb^{-1} = ab^{-1}c,\quad c^6 = 1$.
\end{center}
Therefore any $\sigma\in \mathrm{Aut}^+(L)$ can be written in the form $a^ib^jc^k$ for $i, j \in \mathbb{Z}$ and $0\leq k \leq 5$. This form is unique, so $\mathrm{Aut}^+(L)$ is generated by $a$, $b$ and $c$ subject only to the six relations above.

In fact, $\overline{\Gamma}$ fixes the similar sublattice $R = (\omega + 1)L$, so $\overline{\Gamma}\subset \mathrm{Aut}^+(R)$. In the same way as above, $\mathrm{Aut}^+(R)$ is generated by a translation $A: z\mapsto z+(\omega+1)$, a translation $B: z\mapsto z+(\omega + \omega^2)$, and the rotation $c$, which are subject to the same relations as above on replacing $a$ and $b$ by $A$ and $B$ respectively. By direct computation,
\begin{center}
$\quad \overline{T} = c,\quad \overline{S} = a^{-1}b^2c^3,\quad A = \overline{T}^2\overline{S}\overline{T}^{-2},\quad B = \overline{S}\overline{T}^3$.
\end{center}
So in fact $\overline{\Gamma} = \mathrm{Aut}^+(R)$. Now substitute the relations above into the presentation of $\mathrm{Aut}^+(R)$ to get the desired result.

Let $N$ be the minimal normal subgroup of $\mathrm{SL}_2(\mathbb{Z})$ which contains $T^6$.  The results we will need are summarised in the following corollary.

\cor
\begin{enumerate}
 	\item $\overline{\Gamma} \cong \overline{N}\backslash \Gamma$, where $\overline{N}$ is the image of $N$ under the map $\mathrm{SL}_2(\mathbb{Z}) \to \Gamma$. 
	\item Let $\overline{F}$ be the interior of the triangle with corners $0, \ \omega, \omega^2$. Then $\overline{F}$ is a fundamental domain for the action of $\overline{\Gamma}.$
	\item $ \overline{N} \backslash \mathbb{H} \simeq \B$ as $\Gamma$--sets, where the $\Gamma$--action is given by passing to $\overline{\Gamma}$ via the representation defined earlier and acting on $\B$ as an affine transformation.
\end{enumerate}
\pf These are essentially simple consequences of the proof above. To see that 
$\overline{F}$ is a fundamental domain for the action of $\overline{\Gamma}$ on $\B$ observe that $\overline{T}^k\overline{S}\overline{T}^k$ is a translation by $\omega^k\alpha$.  Pre-multiplying by an appropriate $\overline{T}^h$ then rotates $\overline{F}$ to the required position.

\subsection*{Remark}
The group $N$ is an interesting group in its own right.  As the fundamental group of $\B$, a plane with an infinite number of holes, it is a free group on an infinite number of generators.  Explicit calculation shows that it is generated by matrices of the form $1-6A$ where 

\begin{center}
A = $\begin{pmatrix} -xy&x^2\\-y^2&xy \end{pmatrix}$.
\end{center}
These matrices satisfy $A^2 = 0$, and if $B = \begin{pmatrix} -pq&p^2\\-q^2&pq \end{pmatrix}$, then \\
\[
AB + BA = (xq-py)^2 I.
\]
These matrices then have a symmetric bilinear form.  Multiplication in the algebra generated by these matrices thus behaves like an integral analogue of a Clifford algebra.  As such, its representation theory is interesting in its own right.  

\section{The definition of $\varphi$ }
\subsection{Local definition}

The process here is first to observe that Klein's $J$ invariant maps a fundamental domain for $\Gamma$ bijectively onto the whole of the complex plane, and second to observe that $(\wp')^2$, the square of the derivative of the Weierstrass $\wp$ function similarly sends a fundamental domain for the action of $\overline{\Gamma}$ bijectively onto $\mathbb{C}$.  The results about the Klein's $J$ invariant are standard; the only novelty is observing a parallel structure with respect to its compatibility with $\Gamma$, as in figure 1.  As the results are all classical, we state them with reference to standard literature.

Let $F$ to be the standard fundamental domain for the action of $\Gamma$ on the upper half plane. Let $J$ denote the Klein invariant $J$ map,

\begin{center}
$\displaystyle J(z) = \frac{4}{27}\cdot \frac{(1 - \lambda(z) +\lambda(z)^2)^3}{\lambda(z)^2(1-\lambda(z)^2)}$,
\end{center}
where $\lambda$ is defined by
\begin{center}
$\displaystyle \lambda(z) = \frac{\wp\left(\tfrac{z+1}{2}\right) - \wp\left(\tfrac{z}{2}\right)}{\wp\left(\tfrac{1}{2}\right) - \wp\left(\tfrac{z}{2}\right)}$
\end{center}
relative to the pair of periods $\{1,z\}$.  This function has the required property:

\begin{center}\includegraphics[height=40mm]{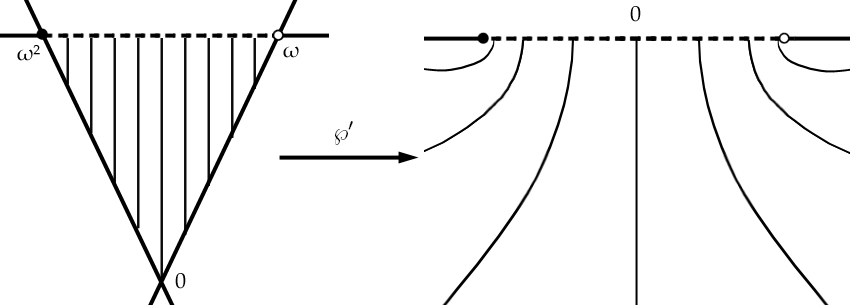}\end{center}

\lem $J$ sends $F$ to $\mathbb{C}$ bijectively and analytically. $J(\omega) = J(\omega^2) = 0$. $J$ is invariant under the action of $\Gamma$.

\pf See Apostol \cite{Apostol:1990wc}, p40.

\lem Let $L = \mathbb{Z}(\omega+1) \oplus \mathbb{Z}(\omega^2-1)$. Then $(\wp'_L)^2$ sends $\overline{F}$ to $\mathbb{C}$ bijectively and analytically. $(\wp'_L)^2$ is invariant under the action of $\overline{\Gamma}$.

\pf See Ahlfors \cite{Ahlfors:1953taa} p89 or McKean and Moll \cite{McKean:1999uf} p268.

Finally, let 
\[
\alpha:\mathbb{C} \to \mathbb{C}, \ \ \alpha(z) = (z-1)/\wp'(\omega)^2.
\]
Putting these together:

\thm \label{junk}
\begin{enumerate}
	\item The map $\varphi|_F$, defined to be the composite
\begin{center}$F \xrightarrow{J} \mathbb{C} \xrightarrow{\alpha} \mathbb{C} \xrightarrow{((\wp'_L)^{2})^{-1}} \overline{F}$,\end{center}
is bijective and analytic.
	\item The extension of $\varphi|_F$ extends to an analytic covering map
	\[
	\varphi: \mathbb{H} \to \B.
	\]
	\item $\varphi$ is a map of $\Gamma$--sets.
\end{enumerate}

\pf By the results quoted above, $\varphi|_F$ is analytic when restricted to $F$. That it is bijective onto $\overline{F}$ does require careful study of where each element of the boundary of $F$ is sent, but all required information is standard knowledge. (See Apostol \cite{Apostol:1990wc} p40ff and McKean and Moll \cite{McKean:1999uf} p89).   

The map $\varphi$ extends to all of $\mathbb{H}$ since $J$ is defined on all of $\mathbb{H}$ and the square root is defined on all of $\mathbb{C} \setminus \{t < 0 | t\in \mathbb{R}\}$ once a branch ($\sqrt{1} = 1$) is chosen.

To see that $\varphi$ is a map of $\Gamma$--sets, the strategy is to check that this is true for $T$ and $S$ on the boundary $\partial F$ of $F$.  By its construction $\varphi(T \partial F) = \overline{T}\varphi(\partial F)$ and $\varphi(S \partial F) = \overline{S}\varphi(\partial F)$. This is not quite enough; we need to check that this is true for all points $u$ in the $\partial F$. But pointwise equality follows since for all $u$ in the boundary,
\[
 (\wp'(\varphi (gu)))^2 = J(gu)  = J(u)  = \wp'(\varphi(u)).
\]
for $g = T, \ S$. The map $(\wp')^2$ is a sixfold covering map; away from zeros and poles it is locally an isomorphism and can be inverted.  Continuity identifies the sheet in which $\varphi(Gu)$ takes its value.

Now let $U$ be an open set containing $F$ and $\partial F$. Then on $gU \bigcap U$, we have that $\varphi(gu) = \overline{g}\varphi(u)$, where $\overline{g} = \psi(g)$.  The function $\varphi(gu) - \overline{g}\varphi(u)$ is analytic on $\mathbb{H}$ and constant (zero) on $U \bigcap \overline{g}\varphi(U)$, hence on all of $\mathbb{H}$.  The result follows since $S$ and $T$ generate $\Gamma$.

\subsection{$\varphi$ and modular functions.} 

The derivative of $\varphi$ has a transformation property similar to that of modular functions.  Observe that

\begin{eqnarray*}
\frac {d}{dz} \varphi (gz) &=& \varphi'(g(z))g'(z) \\
				&=& \frac {d}{dz} \overline{g}(\varphi(z)) \\
				&=& \ell(g)\varphi'(z)
\end{eqnarray*}
where the affine transformation $\overline{g}$ is given in terms of its linear and constant parts:
\[
\overline{g} = \ell(g) + c(g).
\]
But then
\[
\varphi'(g(z)) = \ell(g)\frac{1}{g'(z)}\varphi'(z)
\]
Thus $\varphi'$ transforms as a modular function of weight 2 rotated by a sixth root of unity.  Conventional modular functions then occur as homogeneous polynomials in $(\varphi')^6$, for example.

\section{Group-theoretic description of $\varphi$}

The map $\varphi$ has a simple presentation in terms of double cosets of $\mathrm{SL}_2(\mathbb{R})$.  The choice of a point $p\in \mathbb{H}$ identifies $\mathbb{H}$ with $\mathrm{SL}_2(\mathbb{Z})/\mathrm{Stab}(p)$, where $\mathrm{Stab}(p)$ is the stabiliser of $p$.

Let $N$ be the minimal normal subgroup of $ \mathrm{SL}_2(\mathbb{Z})$ containing $T^6$. Then $N \backslash\mathrm{SL}_2(\mathbb{R})/S_p$ is an $\mathrm{SL}_2(\mathbb{Z})$--space with action given by

\[
gNh\mathrm{Stab}(p) = N(gh)\mathrm{Stab}(p).
\]
The following theorem is a corollary of the proof of theorem \ref{junk}.

\thm $\B \simeq  N \backslash\mathrm{SL}_2(\mathbb{R})/\mathrm{Stab}(p)$  as conformally equivalent $ \mathrm{SL}_2(\mathbb{Z})$--spaces. 

\pf We need to establish that

\[
\B = \{Ng\mathrm{Stab}(p)\ | g \in SL_2(\mathbb{R}) \}.
\]

Let $\pi_p: SL_2(\mathbb{R}) \to \mathbb{H}$ denote the projection $g\mapsto g(p)$, and let $\widetilde{U}$ denote $\pi_p^{-1}(U)$, where $U$ is again an open set containing $F$ and $\partial F$. The final part of the theorem \ref{junk} shows that $\varphi|_U:U \to \overline{U}$ is a covering map. Moreover, $U$ can be chosen such that $T^6U \cap U = \emptyset$.  The consequence of this, together with the fact that $\varphi$ is a map of $\Gamma$--sets, is that $\varphi^{-1}(\overline{U})$ is a disjoint union of the translates $nU$ for $n$ in $N$.  Thus

\[
\overline{U} =  \{Ng\mathrm{Stab}(p)\ | g \in  \widetilde{U}\}.
\]
The proof of the final part of \ref{junk} identifies $\B$ with $\bigcup \overline{g}\overline{U}$ as required.

The final result then can be summarised in the following diagram, where the projection $\Nproj $ sends elements to their right cosets of $N$ in $\Gamma$. (Recall also that, on $\Gamma$, $\Nproj = \psi$.)

\begin{center}
$\xymatrixcolsep{3pc}\xymatrix{
\Gamma\times \mathrm{SL}_2(\mathbb{R})/\mathrm{Stab}(p)\ar[rrr]^-{\text{action}} \ar[ddd]_{1 \times \Nproj} &&& \mathrm{SL}_2(\mathbb{R})/\mathrm{Stab}(p) \ar[ddd]^{\Nproj}\\
&\Gamma\times\mathbb{H}\ar[r]^-{\text{action}} \ar[d]_{\Nproj\times \varphi} \ar@{<->}[ul]^\sim & \mathbb{H}\ar[d]^\varphi \ar@{<->}[ur]_\sim& \\
&\overline{\Gamma}\times \B\ar[r]_-{\text{action}} \ar@<+.3ex>[dl]^{\text{any lift} \times \sim} & \B \ar@{<->}[dr]^\sim& \\
\Gamma\times N\backslash \mathrm{SL}_2(\mathbb{R})/\mathrm{Stab}(p)\ar[rrr]_-{\text{action}} \ar@<+.3ex>[ur]^{\Nproj\times\sim} &&& N\backslash \mathrm{SL}_2(\mathbb{R})/\mathrm{Stab}(p)\\
}$
\end{center}

\subsection*{Remarks}
	This result suggests parallels with the exponential map.  $N$ plays role of $\mathbb{Z}$, $\mathbb{H}$ plays the role of $\mathbb{R}$, and $\B$ the role of the circle. This raises further questions.
	
\begin{enumerate}
	\item A representation of $N$ can be induced to give a representation of $\SL_2(\mathbb{R})$. Given a representation of $N$ - this discrete analogue of an infinite dimensional Clifford algebra - what is the corresponding representation of $ \mathrm{SL}_2(\mathbb{R})$?  How does the representation theory of these two groups correspond?
	\item What is the number--theoretic significance of the equivalence classes of rationals under the action of $N$?
	\item  $ \mathrm{SL}_2(\mathbb{Z})$ acts on all of $\mathbb{C}$.  The equivalence classes of the real numbers under the action of $N$ suggests an interpretation of the geodesics in $\mathbb{H}$ connecting elements within the same equivalence classes as harmonics. How can these be exploited for the purposes of solving differential equations?
\end{enumerate}

\section{Representations of $\mathrm{SL}_3(\mathbb{C})$ and the interpretation of $\wp$ and $\lambda$.}

From the modular function theory point of view, lattices are subgroups $ \lbrace (m\omega_1,n\omega_2)| (m,n) \in \mathbb{Z} \times \mathbb{Z}, \ (\omega_1,\omega_2) \in \mathbb{C}\rbrace$ of $\mathbb{C}$. From the point of view of the representation theory of $\mathfrak{sl}_3(\mathbb{C})$, a diagonal element 
\[
H = \begin{pmatrix}
h_1& 0 & 0 \\
0 & h_2& 0 \\
0 & 0 & h_3 \end{pmatrix}
\]
in the Lie algebra satisfies $h_1+h_2+h_3=0$, and hence setting $\omega_1 = h_1 - h_3,  \ \ \omega_2 = h_2 - h_1$ defines a lattice as in the picture below.

\begin{center}\includegraphics[height=45mm]{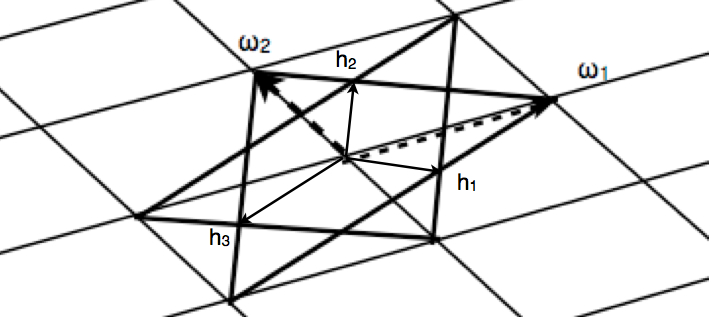}\end{center}

Equivalently, given a pair of generators $(\omega_1,\omega_2)$, define 

\begin {eqnarray*}
H(\omega_1,\omega_2)& =& \begin{pmatrix}
h_1(\omega_1,\omega_2)& 0 & 0 \\
0 & h_2(\omega_1,\omega_2)& 0 \\
0 & 0 & h_3(\omega_1,\omega_2) \end{pmatrix}, \\
	& =& \begin{pmatrix}
\frac{1}{3}\omega_1 - \frac{1}{3}\omega_2& 0 & 0 \\
0 & \frac{1}{3}\omega_1+\frac{2}{3}\omega_2& 0 \\
0 & 0 & -\frac{2}{3}\omega_1-\frac{1}{3}\omega_2 \end{pmatrix}
\end {eqnarray*}

Write $H(\tau)$ and $h_i(\tau)$ for $H(1,\tau)$ and $h_i(1,\tau)$.  For the "standard" lattice, $(\omega_1,\omega_2) = (\omega +1, \omega^2 - 1)$ that we used in establishing the identity $\varphi = (\wp')^2$, the matrix $H = H_{(\omega +1, \omega^2 - 1)}$ is
\[
H = \begin{pmatrix}
1& 0 & 0 \\
0 & \omega^2& 0 \\
0 & 0 & \omega^4 \end{pmatrix}.
\]

For any (finite dimensional) representation $\rho$ of $\mathfrak{sl}_3(\mathbb{C})$, $\rho(H(\tau))$ will be a diagonal operator.  Its trace will then be a function of $\tau$, although for finite dimensional representations this function will not descend to the torus.  The goal of this section is to describe suitable infinite dimensional representations and operators built from $\rho(H(\tau))$ such that the trace of these operators are familiar operators.  Both $\wp$ and modular $\lambda$ can be interpreted in this way.  It may ultimately prove simpler and possibly more illuminating to study these operators than their traces.

\subsection*{Notation for the representation theory of $\mathrm{SL}_3(\mathbb{C})$.}
The diagonal elements in $\mathfrak{sl}_3(\mathbb{C})$ form a commutative sub-Lie algebra $\mathfrak{h}$ called the \emph{Cartan subalgebra} or \emph{maximal torus} of $\mathfrak{sl}_3(\mathbb{C})$. For an element $H \in \mathfrak{h}$, define

\begin{center}
$L_i (H) = L_i \begin{pmatrix}
h_1 & 0 & 0 \\
0 & h_2 & 0 \\
0 & 0 & h_3 \end{pmatrix} = h_i.$
\end{center}

The \emph{weight lattice} is the additive group generated by the $L_i$. The \emph{fundamental weights} are $L_1$ and $-L_3$, and the \emph{simple roots},  $\alpha_1 $, $ \alpha_2$, are given by

\[
\alpha_1 = L_1 - L_2, \ \ \ \alpha_2 = L_2 - L_3.
\]

The simple roots generate the root lattice as a subgroup of the weight lattice. Let $V$ be an irreducible finite-dimensional representation of $\mathfrak{sl}_3(\mathbb{C})$.  The theory of semi-simple Lie algebras guarantees that the diagonal elements of $\mathfrak{sl}_3(\mathbb{C})$ act diagonally on $V$ as well, with the common eigenvalues then appearing as $\mathbb{Z}$ linear combinations of the $L_i$. While no theory guarantees that this convenient state of affairs persists for infinite dimensional representations, some representations do have this property, particularly those derived from polynomials on $\mathbb{C}$.

For any representation $\rho$, the Killing form is the symmetric bilinear invariant form defined by
\[
(X,Y)_{\rho} = \mathrm{tr}\rho(X)\rho(Y)
\]
for elements $X, \ Y$ in $\mathfrak{sl}_3(\mathbb{C})$. This provides an inner product on the real span of the weight lattice.  More generally, leaving matters of convergence aside temporarily, this Killing form is defined for all linear operators on a vector space. Again ignoring matters of convergence, the Killing form is multiplicative with respect to tensor products: if $\rho_1, \ \rho_2$ are two representations,
\[
(*,*)_{\rho_1 \otimes \rho_2} = (*,*)_{\rho_1}(*,*)_{\rho_2}.
\] 

The strategy for the rest of the section is to introduce accessible infinite dimensional representations for  $\mathfrak{sl}_3(\mathbb{C})$.

\subsection{Some useful representations of $\mathfrak{sl}_3$}

Accessible representations of  $\mathfrak{sl}_3(\mathbb{C})$ arise as polynomials on $\mathbb{C}^3$.  Let $\mathbb{C}^3$ have $\{e_1 = x,e_2=y,e_3=z\}$ as basis.  $\mathfrak{sl}_3$ acts on $\mathbb{C}[e_1,e_2,e_3]$ via

\[
E_{ij} f = e_i\partial_jf.
\]

Let $S^n$ denote the homogeneous polynomials of degree $n$.  Those of degree $3n$ will be of particular interest. 

Infinite-dimensional representations in general need not be so neatly arranged, but the space of Laurent polynomials is an example of an infinite-dimensional representation which shares with the finite-dimensional representations the property that $\mathfrak{h}$ acts diagonally and the functionals representing eigenvalues lie in the weight lattice.

Define

\[
W =\left\langle f_{(a,b,c)} := \ \left(\frac {e_1}{e_2}\right)^a\left( \frac {e_2}{e_3}\right)^b \left(\frac{e_3}{e_1}\right)^c  : a, b, c \in \mathbb{Z}_{\geq 0} \right\rangle.
\]

Direct calculation shows that each basis element above spans the one-dimensional weight space $W_{(a,b,c)}$ with weight

\begin{center}
$(a-c)L_1 \ + \ (b-a)L_2 \ + \ (c-b)L_3 $.
\end{center}

Thus the set of weights for this representation is the entire root lattice, with each root occurring with multiplicity one.

This representation can be varied slightly by allowing one of $\{a, \ b, \ c\}$ to take half-integer values. In this case the eigenvalues are represented by linear functionals which no longer lie on the weight lattice.  Define
\[
W(1) = \left\langle f_{(a,b,c)} := \ \left(\frac {e_1}{e_2}\right)^a\left( \frac {e_2}{e_3}\right)^b \left(\frac{e_3}{e_1}\right)^c : a, c \in \mathbb{Z}_{\geq 0}, b \in  \mathbb{Z}_{\geq 0} + \tfrac{1}{2} \right\rangle,
\]
and $W(2), W(3)$ similarly. The weight of $f_{(a,b,c)}$ is again
\[
(a-c)L_1 \ + \ (b-a)L_2 \ + \ (c-b)L_3,
\]
but in $W(i)$, the coefficient of $L_i$ is integral, and that of $L_j$ is half integral for $i \ \neq \ j$.  Write
\[
\rho_0 :  \mathfrak{sl}_3 \to End(W)
\]
and
\[
\rho_i : \mathfrak{sl}_3 \to End(W(i))
\]
for the diagonal action of $\mathfrak{sl}_3$ on $W(i)$.

\subsection{Modular $\lambda$ re-interpreted via the Killing form}

Given any operator $\rho(H)$, its features of interest are its invertibility and its kernel. The operator $(\rho(H) - z)^{-1}$ is therefore a natural choice for study. The Weierstrass $\wp$ function appears to be computing the difference between Killing form of this operator for $\rho = \rho_0$ and the unperturbed operator $\rho_0(H)$ where $H$ is the matrix corresponding to our standard lattice:

\begin{eqnarray}
\nonumber \wp(z, \omega +1, \omega^2-1) &=& tr (\rho_0(H-z)^{-1})^2 - tr(\rho_0(H)^{-1})^2 \\
\nonumber 						&=& ((H-z)^{-1}, (H-z)^{-1})_{\rho_0} - (H^{-1}, H^{-1})_{\rho_0} 
\end{eqnarray}

Of course modular $\lambda$ can then be expressed in terms of the representation $W$, but it there is also the intriguing possibility of representing modular  $\lambda$ in terms of the representations $\rho_i$. This suggests a supersymmetric interpretation.  For simplicity, write $\rho_i(H(\tau)^{-1}) = H_i^{-1}$, and thus write

\[
(H_i^{-1},H_i^{-1}) = (\rho_i(H(\tau)^{-1}),\rho_i(H(\tau)^{-1}))_{\rho_i}.
\]

Then
\[
\lambda(\tau) = \frac{(H_3^{-1},H_3^{-1}) - (H_2^{-1},H_2^{-1})}{(H_1^{-1},H_1^{-1})-(H_2^{-1},H_2^{-1})}.
\]

\subsection{Weierstrass $\wp$ as a limit of traces of finite dimensional representations.}

Since $\wp(z,\omega + 1, \omega^2 - 1)$ is absolutely convergent, the order of summation is unimportant.  Therefore it is possible to choose sequences of finite dimensional representations of $\mathfrak{sl_3}$ which in the limit include the entire root lattice.  For example the representation on polynomials of order $3n$, $\Gamma_{3nL_1}$, of highest weight $3nL_1$, includes all elements of the root lattice inside a triangle with vertices at $\{3n,3n\omega^2,3n\omega^4\}$. Notice that these representations extend to a representation of $\mathfrak{gl}_3$. Then, since we can be confident of convergence, 
\[
\wp(z; \omega + 1,\omega^2 -1) =  \lim_{n\to \infty} ((H-z)^{-1},(H-z)^{-1})_{\Gamma_{3nL_1}} - (H^{-1},H^{-1})_{ \Gamma_{3nL_1}}.
\]
The representations dual to these, $\Gamma_{-3nL_3}$, with highest weight $-3nL_3$, would serve equally well. In matters of computation, it may be preferable to average over these two representations, in order profit from the additional symmetry afforded by duality.  Other sequences of representations, for example those of highest weight $\Gamma_{nL_1 -nL_3}$, can be used by taking the limit of differences.

\bibliographystyle{plain}
\bibliography{bib1.bib}

\end{document}